\documentclass[11pt,a4paper]{amsart}

\usepackage{amssymb,amsmath,amsthm}
\usepackage[all]{xy}

\theoremstyle{plain}
\newtheorem{thm}{Theorem}[section]
\newtheorem{prop}[thm]{Proposition}
\newtheorem{lem}[thm]{Lemma}
\newtheorem{cor}[thm]{Corollary}

\theoremstyle{definition}
\newtheorem{dfn}[thm]{Definition}
\newtheorem{rem}[thm]{Remark}
\newtheorem{exam}[thm]{Example}

\numberwithin{equation}{section}

\makeatletter
\renewenvironment{proof}[1][\proofname]{\par
  \pushQED{\qed}%
  \normalfont \topsep6\p@\@plus6\p@\relax
  \trivlist
  \item[\hskip\labelsep
	\bfseries
    #1\@addpunct{.}]\ignorespaces
}{%
  \popQED\endtrivlist\@endpefalse
}
\makeatother

\DeclareMathOperator{\Ext}{Ext}
\DeclareMathOperator{\End}{End}
\DeclareMathOperator{\wt}{wt}
\DeclareMathOperator{\rad}{rad}
\DeclareMathOperator{\soc}{soc}
\DeclareMathOperator{\ch}{ch}
\DeclareMathOperator{\gch}{gch}

\newcommand{\z}{\mathbb{Z}}
\newcommand{\C}{\mathbb{C}}

\newcommand{\M}{\mathfrak{M}}
\newcommand{\La}{\mathfrak{L}}

\newcommand{\g}{\mathfrak{g}}
\newcommand{\h}{\mathfrak{h}}
\newcommand{\bo}{\mathfrak{b}}
\newcommand{\n}{\mathfrak{n}}

\newcommand{\bfi}{\mathbf{i}}
\newcommand{\cl}{\mathrm{cl}}
\newcommand{\dext}{D_{\mathbf{i}}^{\mathrm{ext}}}

\newcommand{\e}[1]{\tilde{e}_{#1}}

\newcommand{\gauss}[2]{\genfrac{[}{]}{0pt}{0}{#1}{#2}}

\title[Loewy series of Weyl modules and quiver varieties]{Loewy series of Weyl modules and the Poincar\'{e} polynomials of quiver varieties}

\author{Ryosuke Kodera}
\author{Katsuyuki Naoi}

\address[Ryosuke Kodera]{Graduate School of Mathematical Sciences, The University of Tokyo, 3-8-1 Komaba, Meguro-ku, Tokyo 153-8914, Japan}
\email{kryosuke@ms.u-tokyo.ac.jp}

\address[Katsuyuki Naoi]{Graduate School of Mathematical Sciences, The University of Tokyo, 3-8-1 Komaba, Meguro-ku, Tokyo 153-8914, Japan}
\email{naoik@ms.u-tokyo.ac.jp}


\begin{document}

\begin{abstract}
We prove that a Weyl module for the current Lie algebra associated with a simple Lie algebra of type $ADE$ is rigid, that is, it has a unique Loewy series.
Further we use this result to prove that the grading on a Weyl module defined by the degree of currents coincides with another grading which comes from the degree of the homology group of the quiver variety.
As a corollary we obtain a formula for the Poincar\'{e} polynomials of quiver varieties of type $ADE$ in terms of the energy functions defined on the crystals for tensor products of level-zero fundamental representations of the corresponding quantum affine algebras. 
\end{abstract}

\maketitle

\section{Introduction}

The purpose of this article is to study the graded structures of Weyl modules for the current Lie algebra associated with a simple Lie algebra of type $ADE$ and its applications.

Weyl modules are defined for a current Lie algebra, a loop Lie algebra, or more generally, a generalized current Lie algebra and play an important role in the study of finite-dimensional modules over these Lie algebras.
They were originally introduced by Chari and Pressley in \cite{MR1850556} for a loop Lie algebra with an application in their mind to the representation theory of the corresponding quantum loop algebra.
This notion was extended by Feigin and Loktev in \cite{MR2102326} to the case of a generalized current Lie algebra (the current Lie algebra on an affine variety, they called), including an untwisted multivariable current (or loop) Lie algebra.
Applications of Weyl modules have been made by Chari and Moura in \cite{MR2078944} to the block decomposition of the category of finite-dimensional modules over a loop Lie algebra, and by the first author in \cite{MR2657446} to a problem analogous to \cite{MR2078944} for a generalized current Lie algebra, for instance.

When we focus on the current Lie algebra associated with a simple Lie algebra $\g$, the Weyl module $W(\lambda)$ for a dominant integral weight $\lambda \in P_+$ of $\g$ admits one remarkable property; it has a \emph{graded} module structure.
This grading helps us to investigate a detailed structure of the Weyl module.
Moreover it has gradually turned out that the gradings of Weyl modules are connected with some interesting objects.
One connection we emphasize here is that with the one-dimensional sums associated with certain crystals, where such a connection has been found in a series of recent works as \cite{MR2271991} by Chari and Loktev, \cite{MR2323538} by Fourier and Littelmann, and \cite{naoi} by the second author.
As another one, from results of Ardonne and Kedem in \cite{MR2290922} and Di Francesco and Kedem in \cite{MR2428305}, we deduce that the gradings are also related to the fermionic forms.
This observation solves the $X=M$ conjecture for particular cases (See \cite[Section~9]{naoi} for a detailed explanation).
These results suggest that the graded structure of a Weyl module itself is interesting and worth studying. 

In this article we study in more detail Weyl modules for the current Lie algebra associated with a simple Lie algebra $\g$ of type $ADE$.
For the case of type $ADE$, the Weyl module $W(\lambda)$ is isomorphic to two other important modules.
\begin{itemize}
\item The Weyl module $W(\lambda)$ is isomorphic to a Demazure module for the affine Lie algebra $\hat{\g}$ (\cite{MR2271991}, \cite{MR2323538}).

\item The Weyl module $W(\lambda)$ is isomorphic to the standard module $M(\lambda)$ defined as the homology group of the Lagrangian quiver variety $\La(\lambda)$ (Proposition~\ref{prop:standard}).
\end{itemize}
Thus it is expected that Weyl modules for type $ADE$ have more specific properties than those for other types and that the identification with the other modules above is useful for understanding of them.  
This expectation will turn out to be true.

Our first main result is to determine the Loewy structure of a Weyl module.
Recall that a Loewy series of a module of finite length is by definition a semisimple filtration which has the smallest length.
It is a fundamental problem to study the Loewy structure of a module, especially to determine two standard Loewy series: the radical series and the socle series.
We say that a module is rigid if its radical series and socle series coincide, which implies that its Loewy series is unique.
Because of this property, rigid modules are fairly easy to understand.
We prove that Weyl modules for type $ADE$ are rigid as an application of the identification with the Demazure modules together with the graded structures.  

\begin{thm}[Proposition~\ref{prop:radical}, Proposition~\ref{prop:socle} and Theorem~\ref{thm:rigid}]\label{thm:rigid1}
For the Weyl module $W(\lambda)$, the grading filtration, the radical series and the socle series coincide.
In particular, $W(\lambda)$ is rigid.
\end{thm}

We remark that a Weyl module for type $BCFG$ is not rigid in general (See Example~\ref{exam:couterexample}).
Thus rigidity is a specific phenomenon only for type $ADE$.

As well as the Weyl module $W(\lambda)$, the standard module $M(\lambda)$ has a graded structure which comes from the degree of the homology group of the quiver variety.
The second main result is to prove that their gradings coincide.

\begin{thm}[Theorem~\ref{thm:coincidence}]\label{thm:coincidence1}
Under the isomorphism $M(\lambda) \cong W(\lambda)$, the gradings on the both sides coincide.
\end{thm}

It should be noticed that Theorem~\ref{thm:coincidence1} seems to be known but has not been in the literature. 
See Remark~\ref{rem:Varagnolo} for a more precise explanation.
In this article, we give a proof of this fact by the representation theory, using rigidity of the Weyl module stated in Theorem~\ref{thm:rigid1}.

Let us mention two corollaries of Theorem~\ref{thm:coincidence1}.
We see by Theorem~\ref{thm:coincidence1} that the graded characters of $M(\lambda)$ and $W(\lambda)$ are equal.
The former is nothing but the generating function of the Poincar\'{e} polynomials of the quiver varieties, which follows immediately from the definition of the standard module and the grading on it.
The latter is expressed by the second author in \cite{naoi} in terms of the degree function $D$ defined on the crystal $\mathbb{B}(\lambda)_{\mathrm{cl}}$ of Lakshmibai-Seshadri paths of shape $\lambda$ modulo imaginary roots, which are studied by Naito and Sagaki in \cite{MR2146858} and \cite{MR2474320}, or equivalently by using the normalized one-dimensional sums $X(\lambda, \mu, t)$ for dominant integral weights $\lambda, \mu \in P_+$.  
As a consequence we obtain the following formula, where $\La(\alpha, \lambda)$ denotes the Lagrangian quiver variety associated with $\lambda \in P_+$ and $\alpha \in Q_+$, and $d_{\alpha}$ stands for the twice of the dimension of $\La(\alpha, \lambda)$.

\begin{cor}[Corollary~\ref{cor:Poincare}]
We have
\begin{align*}
\sum_{\alpha \in Q_+} \sum_{k=0}^{d_{\alpha}} \dim H_k(\La(\alpha, \lambda))t^{d_{\alpha}-k} e^{\lambda - \alpha} &= \sum_{b \in \mathbb{B}(\lambda)_{\cl}}t^{-2D(b)}e^{\wt b}\\
&= \sum_{\mu \in P_+} X(\lambda, \mu, t^{-2})\ch V(\mu).
\end{align*}
In particular, we have
\[
\sum_{k=0}^{d_{\alpha}} \dim H_k(\La(\alpha, \lambda))t^{d_{\alpha}-k}= \sum_{\substack{b \in \mathbb{B}(\lambda)_{\cl}\\ \wt b = \lambda - \alpha}}t^{-2D(b)}.
\]
\end{cor}

We note that some formulas for the Poincar\'{e} polynomials of quiver varieties have been established (\cite{MR2651380}, \cite{moz}, \cite{MR2144973} for instance).
However those in terms of crystals have not been in the literature as far as the authors know.

We also obtain the equality between the Kazhdan-Lusztig type polynomial $Z_{\lambda \mu}(t)$ for the quiver variety and the one-dimensional sum.

\begin{cor}[Corollary~\ref{cor:KL}]
We have
\[
	Z_{\lambda \mu}(t) = X(\lambda, \mu, t^{-2}).
\]
\end{cor}

The article is organized as follows.
In Section~\ref{section:Loewy} we recall basic facts on Loewy series of modules.
We also provide a key lemma to the proof of rigidity of Weyl modules.
In Section~\ref{section:Weyl} we give the definition of Weyl modules and prove that they are rigid.
In addition we determine the socle and the Loewy length of a Weyl module.
The most part of Section~\ref{section:quiver} and \ref{section:cohomological} is a summary of known results.
After recalling properties of quiver varieties in the first half of Section~\ref{section:quiver}, we define standard modules in the latter half.
In Section~\ref{section:cohomological} we introduce a grading on a standard module and interpret it in terms of sheaves on a quiver variety.
In Section~\ref{section:coincidence} we prove that the gradings on a Weyl module and a standard module coincide and deduce its corollaries.
This section also contains discussions on some related subjects such as Lusztig's fermionic conjecture, the $X=M$ conjecture and Nakajima's result on quiver varieties of type $A$. 
 
\subsection*{Acknowledgments}

The authors would like to thank Hiraku Nakajima for his valuable suggestions, especially pointing out that Theorem~\ref{thm:coincidence1} is proved more directly.
They also thank Yoshihisa Saito who read the manuscript carefully and gave many helpful comments.
The both authors are supported by Research Fellowships of the Japan Society for the Promotion of Science for Young Scientists.

\section{Loewy series and the grading filtrations}\label{section:Loewy}

\subsection{Loewy series of modules}

We recall basic notion such as radical, socle and Loewy series of modules and their properties in this subsection.

Let $A$ be a ring and $M$ an $A$-module of finite length.
A filtration of $A$-modules on $M$ is said to be semisimple if its each successive quotient is semisimple.
The radical of $M$, which is denoted by $\rad M$, is the smallest submodule of $M$ such that the quotient is semisimple.
We put $\rad^0 M = M$ and for $k \geq 1$, define $\rad^k M$ inductively by $\rad^k M = \rad(\rad^{k-1} M)$.
This defines a semisimple filtration on $M$ called the \emph{radical series}.
The socle of $M$, which is denoted by $\soc M$, is the largest semisimple submodule of $M$.
We put $\soc^0 M = 0$ and for $k \geq 1$, define $\soc^k M$ inductively so that $\soc (M/\soc^{k-1} M) = \soc^{k} M/\soc^{k-1} M$.
This defines a semisimple filtration on $M$ called the \emph{socle series}.
It is easy to show that for any semisimple filtration
\[
	0 = F^lM \subseteq F^{l-1}M \subseteq  \cdots \subseteq F^{1}M \subseteq F^{0}M = M
\]
on $M$,
\[
	\rad^{k} M \subseteq F^{k} M \subseteq \soc^{l-k}M
\]
holds for each $k$.
This immediately implies that the lengths of the radical series and the socle series are equal and that the length of any semisimple filtration on $M$ is greater than or equal to it.  A \emph{Loewy series} of an $A$-module $M$ of finite length is defined to be a semisimple filtration on $M$ which has this smallest length.
The radical series and the socle series are Loewy series by the definition.
The length of each Loewy series of $M$ is called the \emph{Loewy length} of $M$.

\begin{dfn}
Let $A$ be a ring and $M$ an $A$-module of finite length.
We say that $M$ is \emph{rigid} if its radical series and socle series coincide.\end{dfn}

The following is obvious from the above argument.	

\begin{prop}
Let $A$ be a ring and $M$ an $A$-module of finite length.
If $M$ is rigid then $M$ has a unique Loewy series.
\end{prop}

\subsection{Filtrations on graded modules}

In this article, by grading we always mean a $\z$-grading and by a positively graded ring or module we mean that it is graded by nonnegative integers.

Let $A$ be a graded ring.
If $A$ is positively graded then a graded $A$-module $M$ is endowed with the filtration defined by $F^kM = \bigoplus_{s \geq k} M_s$.
We call it the \emph{grading filtration}.
Moreover assume that $A$ is a positively graded $\C$-algebra such that any finite-dimensional $A_0$-module is semisimple, and that $M$ is a finite-dimensional graded $A$-module.
Then in this situation, the grading filtration on $M$ is semisimple.

We slightly modify a proposition in \cite[Proposition~2.4.1]{MR1322847} for our setting and state it as follows.
The proof goes without any change.

\begin{lem}\label{lem:BGS}
Let $A$ be a positively graded $\C$-algebra and suppose that any finite-dimensional $A_0$-module is semisimple and that $A$ is generated by $A_1$ as an $A_0$-algebra.
Then for a finite-dimensional graded $A$-module $M$ the following hold.
\begin{enumerate}
\item If $M/\rad M$ is simple then the radical series of M coincides with the grading filtration.
\item If $\soc M$ is simple then the socle series of M coincides with the grading filtration.
\end{enumerate}
\end{lem}

\section{Loewy series of Weyl modules}\label{section:Weyl}

\subsection{Weyl modules for a current Lie algebra}\label{subsection:Weyl}

Let $\mathfrak{a}$ be an arbitrary Lie algebra over $\C$.
The loop Lie algebra associated with $\mathfrak{a}$ is the tensor product $\mathfrak{a} \otimes \C[z,z^{-1}]$ equipped with the Lie algebra structure by $[x \otimes f, y \otimes g] = [x,y] \otimes fg$ for $x,y \in \mathfrak{a}$ and $f,g \in \C[z,z^{-1}]$.
We denote by $\mathfrak{a}[z,z^{-1}]$ this Lie algebra.
The \emph{current Lie algebra} is the Lie subalgebra $\mathfrak{a} \otimes \C[z]$ of $\mathfrak{a}[z,z^{-1}]$ and denoted by $\mathfrak{a}[z]$.
We denote by $z\mathfrak{a}[z]$ the Lie subalgebra $\mathfrak{a} \otimes z\C[z]$.
We regard $\mathfrak{a}$ as a Lie subalgebra of $\mathfrak{a}[z]$ by identifying it with $\mathfrak{a} \otimes 1$.

Let $\g$ be a simple Lie algebra over $\C$.
In this subsection, we impose no assumption on $\g$ while we will assume that $\g$ is of type $ADE$ from the next subsection.
We fix a Cartan subalgebra $\h$ of $\g$ and a Borel subalgebra $\bo$ containing $\h$. 
The nilpotent radical of $\bo$ is denoted by $\n$.
Let $I$ be the index set of simple roots.
The simple roots are denoted by $\alpha_i$ ($i \in I$) and the fundamental weights by $\varpi_i$ ($i \in I$).
We choose Chevalley generators $e_i,h_i,f_i$ ($i \in I$) of $\g$.

Let $P$ be the weight lattice and $P_+$ the set of all dominant integral weights.
Let $Q$ be the root lattice and $Q_+$ its subset consisting of all elements expressed as sums of simple roots with nonnegative coefficients.
For $\lambda, \mu \in P$ we say that $\lambda \geq \mu$ if $\lambda - \mu \in Q_+$. 

The Weyl group of $\g$ is denoted by $W$ and its longest element by $w_0$.
Let $(\, ,\,)$ be the $W$-invariant nondegenerate symmetric bilinear form on $P$ normalized by $(\alpha, \alpha) = 2$ for a long root $\alpha$.

The universal enveloping algebra $U(\g[z])$ of the current Lie algebra $\g[z]$ is graded by the degree of $z$.
It is obvious that $U(\g[z])_0 = U(\g)$.
Note that $U(\g[z])$ satisfies the assumptions in Lemma~\ref{lem:BGS}, namely it is a positively graded $\C$-algebra such that any finite-dimensional $U(\g[z])_0$-module is semisimple and $U(\g[z])$ is generated by $U(\g[z])_1$ as a $U(\g[z])_0$-algebra.

For $\lambda \in P_+$, let $V(\lambda)$ be the finite-dimensional simple $U(\g)$-module with highest weight $\lambda$.
We define the $U(\g[z])$-module structure on $V(\lambda)$ through the projection from $U(\g[z])$ to $U(\g[z])_0 = U(\g)$ and denote it by the same symbol $V(\lambda)$.
It is simple as a $U(\g[z])$-module.

Now we give the definition of Weyl modules for the current Lie algebra $\g[z]$ following a formulation by Chari and Loktev in \cite{MR2271991}.

\begin{dfn}
Let $\lambda$ be an element of $P_+$.
The \emph{Weyl module} $W(\lambda)$ is the $U(\g[z])$-module generated by a nonzero element $v_{\lambda}$ with the following defining relations:
\[
	\n[z]v_{\lambda}=0,
\]
\[
	z\h[z]v_{\lambda}=0,
\]
\[
	hv_{\lambda} = \langle h, \lambda \rangle v_{\lambda}
\]
for $h \in \h$,
\[
	f_i^{\langle h_i, \lambda \rangle + 1}v_{\lambda} = 0
\]
for $i \in I$.
\end{dfn}

Proposition~\ref{prop:top} below is proved by a standard argument and Proposition~\ref{prop:graded} is immediate from the definition.

\begin{prop}\label{prop:top}
The Weyl module $W(\lambda)$ has a unique simple quotient $V(\lambda)$.
\end{prop}

\begin{prop}\label{prop:graded}
The Weyl module $W(\lambda)$ is a graded $U(\g[z])$-module.
\end{prop}

We set a grading on $W(\lambda)$ so that the degree of $v_{\lambda}$ is zero.
Then it is positively graded and we have $W(\lambda)_0 \cong V(\lambda)$ as $U(\g)$-modules.

As stated in \cite[1.2.2~Theorem]{MR2271991}, the following fundamental result  is proved in a way similar to \cite[Theorem~1~(ii)]{MR1850556} for the case of a loop Lie algebra (See also \cite[Theorem~1]{MR2102326} for a more general setting).

\begin{thm}[]\label{thm:finite-dimensional}
The Weyl module $W(\lambda)$ is finite-dimensional.
\end{thm}

We have an immediate consequence of Proposition~\ref{prop:top}, Proposition~\ref{prop:graded} and Theorem~\ref{thm:finite-dimensional} together with Lemma~\ref{lem:BGS} (i).

\begin{prop}\label{prop:radical}
The radical series of $W(\lambda)$ coincides with the grading filtration.
\end{prop}

\subsection{Weyl modules and Demazure modules}

Assume that $\g$ is of type $ADE$ in this subsection.
We review a relation between Weyl modules and Demazure modules when $\g$ is of type $ADE$.

Let $\hat{\g}$ be the untwisted affine Lie algebra associated with $\g$, namely  $\hat{\g} = \g[z,z^{-1}] \oplus \C c \oplus \C d$ as a $\C$-vector space, where $c$ is the canonical central element and $d$ is the degree operator.
The Cartan subalgebra $\hat{\h}$ and the Borel subalgebra $\hat{\bo}$ are given by $\hat{\h} = \h \oplus \C c \oplus \C d$ and $\hat{\bo} = \bo \oplus z\g[z] \oplus \C c \oplus \C d$.

Let $\hat{I} = I \sqcup \{0\}$ be the index set of simple roots of $\hat{\g}$.
Let $\hat{P}$ be the weight lattice and $\hat{P}_+$ the set of all dominant integral weights.
The fundamental weights are denoted by $\Lambda_i$ ($i \in \hat{I}$) and the generator of the imaginary roots by $\delta$.
We regard $P$ as a subset of $\hat{P}$ by $\varpi_i = \Lambda_i - \langle c, \Lambda_i \rangle \Lambda_0$ for $i \in I$.
Then any $\Lambda \in \hat{P}$ is uniquely expressed as $\Lambda = \lambda + l\Lambda_0 + m\delta$ for some $\lambda \in P$ and integers $l, m$.
The integer $l$, which is equal to $\langle c, \Lambda \rangle$, is called the level of $\Lambda$.

The Weyl group of $\hat{\g}$ is denoted by $\hat{W}$.
The bilinear form on $P$ is extended to the $\hat{W}$-invariant nondegenerate symmetric bilinear form $(\, ,\, )$ on $\hat{P}$ satisfying $(\alpha_i,\Lambda_j) = \delta_{ij}$ for $i,j \in \hat{I}$ and $(\Lambda_0, \Lambda_0) = 0$.
The translation $t_{\alpha}$ by $\alpha \in Q$ on $\hat{P}$ is defined by $t_{\alpha}(\Lambda) = \Lambda + \langle c, \Lambda \rangle\alpha - \{(\Lambda, \alpha) + (1/2)(\alpha, \alpha)\langle c, \Lambda \rangle\}\delta$ for $\Lambda \in \hat{P}$.
Since $t_{\alpha} \in \hat{W}$ for any $\alpha \in Q$, the root lattice $Q$ of $\g$ is regarded as a subgroup of $\hat{W}$.
Then $\hat{W}$ is isomorphic to the semidirect product $W \ltimes Q$.

Let $L(\Lambda)$ be the integrable simple $U(\hat{\g})$-module with highest weight $\Lambda \in \hat{P}_+$.
For each $w \in \hat{W}$ the extremal weight space $L(\Lambda)_{w\Lambda}$ is one-dimensional.
 
\begin{dfn}
Let $\Lambda$ be an element of $\hat{P}_+$ and $w$ an element of $\hat{W}$.
The \emph{Demazure module} $L_w(\Lambda)$ is the $U(\hat{\bo})$-submodule of $L(\Lambda)$ generated by the extremal weight space $L(\Lambda)_{w\Lambda}$.
\end{dfn}

In the sequel we consider the Demazure modules which have the $U(\g[z])$-module structures.
For $\lambda \in P_+$ we can take $w \in \hat{W}$ and $\Lambda \in \hat{P}_+$ so that $w\Lambda = w_0\lambda + \Lambda_0$.
The choice of $\Lambda$ is unique.
In this situation, it is known that the Demazure module $L_w(\Lambda)$ is $U(\g[z])$-stable (See \cite[a sentence below Remark~3.5]{naoi} for example).
This $L_w(\Lambda)$ is positively graded by the eigenvalues for the action of the degree operator $d$.
Note that the restriction of $w_0w\Lambda = \lambda + \Lambda_0$ to $\h$ is equal to $\lambda$.

The following theorem, which asserts that a Weyl module is isomorphic to a Demazure module, was proved by Chari and Loktev in \cite[1.5.1~Corollary]{MR2271991} for type $A$ and by Fourier and Littelmann in \cite[Theorem~7]{MR2323538} for type $ADE$.

\begin{thm}[]\label{thm:FL}
Let $\lambda$ be an element of $P_+$.
Take $w \in \hat{W}$ and $\Lambda \in \hat{P}_+$ so that $w\Lambda = w_0\lambda + \Lambda_0$.
Then there exists an isomorphism $W(\lambda) \cong L_w(\Lambda)$ of graded $U(\g[z])$-modules which sends $v_\lambda$ to an extremal weight vector $u_{w_0w\Lambda}$ of weight $w_0w\Lambda$.
\end{thm}

\subsection{Rigidity}

In this subsection we prove that Weyl modules are rigid for the current Lie algebra associated with a simple Lie algebra $\g$ of type $ADE$.
Since we proved in Subsection~\ref{subsection:Weyl} that the radical series of a Weyl module coincides with the grading filtration, it suffices to show that the socle series also coincides with it.
In addition we determine the socle and the Loewy length of a Weyl module.
We assume that $\g$ is of type $ADE$ unless otherwise specified in this subsection.

\begin{prop}\label{prop:simple_socle}
The socle of $W(\lambda)$ is simple.
\end{prop}
\begin{proof}
By Theorem~\ref{thm:FL}, $W(\lambda)$ is isomorphic to the Demazure module $L_w(\Lambda)$ for some $w \in \hat{W}$ and $\Lambda \in \hat{P}_+$.
Since any nonzero $U(\g[z])$-submodule of $L_w(\Lambda)$ contains the highest weight space $L(\Lambda)_{\Lambda}$, the assertion follows.
\end{proof}

The desired property of the socle series follows from the above proposition, Proposition~\ref{prop:graded} and Theorem~\ref{thm:finite-dimensional} together with Lemma~\ref{lem:BGS} (ii).

\begin{prop}\label{prop:socle}
The socle series of $W(\lambda)$ coincides with the grading filtration.
\end{prop}

By Proposition~\ref{prop:radical} and Proposition~\ref{prop:socle} we obtain the following, which is the first main theorem of this article.

\begin{thm}\label{thm:rigid}
The Weyl module $W(\lambda)$ is rigid.
\end{thm}

\begin{rem}
If $\g$ is not of type $ADE$, namely is of type $BCFG$, then the socle series of the Weyl module $W(\lambda)$ does not coincide with the grading filtration in general, while the radical series always coincides with it as proved in Proposition~\ref{prop:radical}.
Thus it possibly has the nonsimple socle and is not rigid.
An example is given in the following.    
\end{rem}

\begin{exam}\label{exam:couterexample}
Let $\g$ be of type $C_2$. The index set $I=\{1,2\}$ is numbered so that $\alpha_1$ is a short root and $\alpha_2$ is long.
Let $\lambda = 2\varpi_1 + \varpi_2$ and $\mu = 2\varpi_2$.
Take $w, w' \in \hat{W}$ and $\Lambda, \Lambda' \in \hat{P}_+$ so that $w\Lambda = w_0\lambda + \Lambda_0$ and $w'\Lambda' = w_0\mu + \Lambda_0$.
By \cite[Theorem~9.3]{naoi}, we see that there exists an exact sequence
\[
	0 \to L_{w'}(\Lambda') \to W(\lambda) \to L_{w}(\Lambda) \to 0
\]
of $\g[z]$-modules.
The following figure expresses the grading on $W(\lambda)$.
\[
\begin{array}{ccccc}
\text{degree} && L_{w}(\Lambda) && L_{w'}(\Lambda') \\ \hline
0 && V(2\varpi_1 + \varpi_2) && \\ \hline
1 && V(2\varpi_1) \oplus V(\varpi_2) && V(2\varpi_2) \\ \hline
2 && V(\varpi_2) && V(2\varpi_1) \\ \hline
3 && && V(0)
\end{array}
\]
Take notice of the composition factors $V(0)$ in degree $3$ and $V(\varpi_2)$ in degree $2$.
It is known that no nontrivial extension between $V(0)$ and $V(\varpi_2)$ exists (See \cite[Proposition~3.1 and Remark~3.5]{MR2657446} for example.
The highest root $\theta$ is given by $2\varpi_1$ here).
Therefore the socle of $W(\lambda)$ contains $V(0) \oplus V(\varpi_2)$ hence is not simple (In this case it is easy to check that the socle coincides with $V(0) \oplus V(\varpi_2)$ exactly).
We also see that $W(\lambda)$ is not rigid.

In this example, the Weyl module $W(\lambda)$ is \emph{not} isomorphic to any Demazure module but has a filtration such that its each successive quotient is isomorphic to a Demazure module.
This phenomenon is common for a general simple Lie algebra as proved by the second author in \cite{naoi}. 
\end{exam} 

At the end of this section we determine the socle and the Loewy length of $W(\lambda)$.

\begin{lem}\label{lem:Demazure}
Let $\lambda$ be an element of $P_+$.
\begin{enumerate}
\item There exists a unique minimal element $\lambda_{\mathrm{min}}$ in $\{\mu \in P_+ \mid \mu \leq \lambda\}$.

\item The element $\lambda_{\mathrm{min}}$ is equal to $\varpi_i$ for some $i \in I$ such that $\langle c, \Lambda_i \rangle = 1$ or is equal to $0$. 

\item Take $w \in \hat{W}$ and $\Lambda \in \hat{P}_+$ so that $w\Lambda = w_0\lambda + \Lambda_0$.
Then $\Lambda$ is equal to $\lambda_{\mathrm{min}} + \Lambda_0 + (1/2)\{(\lambda,\lambda)-(\lambda_{\mathrm{min}}, \lambda_{\mathrm{min}})\}\delta$.
\end{enumerate}
\end{lem}
\begin{proof}
The assertions of (i) and (ii) are well known and are proved by standard arguments.

We prove (iii).
Put $\xi = \lambda - \lambda_{\mathrm{min}} \in Q_+$.
We calculate $t_{-\xi}w_0w \Lambda = t_{-\xi}(\lambda + \Lambda_0)$ by a formula in \cite[(6.5.3)]{MR1104219} as
\begin{align*}
	t_{-\xi}w_0w \Lambda &=t_{-\xi}(\lambda + \Lambda_0)\\
	&= \lambda + \Lambda_0 - \xi + (1/2)\{(\lambda + \Lambda_0, \lambda + \Lambda_0)-(\lambda_{\mathrm{min}}, \lambda_{\mathrm{min}})\}\delta \\
	&= \lambda_{\mathrm{min}} + \Lambda_0 + (1/2)\{(\lambda, \lambda)-(\lambda_{\mathrm{min}}, \lambda_{\mathrm{min}})\}\delta,
\end{align*}
where the third equality holds since we have $(\lambda, \Lambda_0) = (\Lambda_0, \Lambda_0) = 0$.
By (ii), $\lambda_{\mathrm{min}} + \Lambda_0$ is equal to some fundamental weight $\Lambda_i$ and hence $t_{-\xi}w_0w \Lambda$ belongs to $\hat{P}_+$.
Since $\Lambda \in \hat{P}_+$, this implies that $t_{-\xi}w_0w \Lambda = \Lambda$.
The assertion is proved. 
\end{proof}

We use the symbol $\lambda_{\mathrm{min}}$ for the unique minimal element in the sequel.

\begin{prop}\label{prop:length}
Let $\lambda$ be an element of $P_+$.
\begin{enumerate}
\item The socle of $W(\lambda)$ is isomorphic to $V(\lambda_{\mathrm{min}})$.
\item The Loewy length of $W(\lambda)$ is equal to $(1/2)\{(\lambda,\lambda)-(\lambda_{\mathrm{min}}, \lambda_{\mathrm{min}})\} + 1$.
\end{enumerate}
\end{prop}
\begin{proof}
Take $w \in \hat{W}$ and $\Lambda \in \hat{P}_+$ so that $w\Lambda = w_0\lambda + \Lambda_0$.
Then we have $W(\lambda) \cong L_w(\Lambda)$ by Theorem~\ref{thm:FL}.

We prove (i).
Recall that the socle of $W(\lambda)$ is simple by Proposition~\ref{prop:simple_socle}.
Under the isomorphism $W(\lambda) \cong L_{w}(\Lambda)$, the socle of $W(\lambda)$ is isomorphic to the $U(\g[z])$-submodule generated by the highest weight space $L(\Lambda)_{\Lambda}$.
It is isomorphic to $V(\Lambda |_{\h})=V(\lambda_{\mathrm{min}})$.

We prove (ii).
The grading filtration on $W(\lambda)$ gives its Loewy series by Proposition~\ref{prop:radical}.
Hence the Loewy length is equal to $\max\{k \mid W(\lambda)_k \neq 0\} + 1$.
By Lemma~\ref{lem:Demazure} (iii), $\max\{k \mid W(\lambda)_k \neq 0\}$ is equal to $(1/2)\{(\lambda,\lambda)-(\lambda_{\mathrm{min}}, \lambda_{\mathrm{min}})\}$.
\end{proof}

\begin{rem}
The integer $(\lambda,\lambda)-(\lambda_{\mathrm{min}}, \lambda_{\mathrm{min}})$ appearing in the Loewy length of $W(\lambda)$ will be interpreted as the dimension of the nonsingular quiver variety $\M(\lambda-\lambda_{\mathrm{min}}, \lambda)$ and also of the affine quiver variety $\M_0(\lambda)$ in Section~\ref{section:coincidence}.
\end{rem}

\section{Quiver varieties and standard modules}\label{section:quiver}

In the remaining of this article we assume that $\g$ is of type $ADE$.
Any algebraic variety, sometimes variety for short, is assumed to be over $\C$ and not necessarily to be irreducible or connected in the sequel.
The dimension of an algebraic variety always means the complex dimension.

\subsection{Quiver varieties}

\emph{Quiver varieties} were introduced by Nakajima in \cite{MR1302318} attached to graphs with some additional data in terms of hyper-K\"{a}hler quotients.
Later in \cite{MR1604167} they were reformulated by the geometric invariant theory.
We denote by $\M(\alpha, \lambda)$ and $\M_0(\alpha, \lambda)$ the quiver varieties associated with $\lambda \in P_+$ and $\alpha \in Q_+$, where they correspond to $\M(\mathbf{v}, \mathbf{w})$ and $\M_0(\mathbf{v}, \mathbf{w})$ defined in \cite[Section~3]{MR1604167} respectively, by identifying $\alpha = \sum_{i \in I} v_{i} \alpha_i \in Q_+$ with $\mathbf{v} = (v_i)_{i \in I} \in (\z_{\geq0})^I$ and $\lambda = \sum_{i \in I} w_{i} \varpi_i \in P_+$ with $\mathbf{w} = (w_i)_{i \in I} \in (\z_{\geq0})^I$.
Remark that we consider the quiver varieties only for the Dynkin diagram of the simple Lie algebra $\g$ of type $ADE$.

Let us gather basic properties of quiver varieties.
See \cite{MR1302318}, \cite{MR1604167} and \cite{MR1808477} for proofs.
Note that some of the properties stated below hold only for the case of type $ADE$, not for an arbitrary graph.

\begin{itemize}
\item For $\lambda \in P_+$ and $\alpha \in Q_+$, $\M(\alpha, \lambda)$ is a possibly empty nonsingular quasi-projective variety and $\M_0(\alpha, \lambda)$ is an affine variety.

\item There exists a projective morphism $\pi\colon \M(\alpha, \lambda) \to \M_0(\alpha, \lambda)$.

\item The variety $\M(\alpha, \lambda)$ has a symplectic structure.
\end{itemize}

The variety $\M_0(\alpha, \lambda)$ has a distinguished point denoted by $0$.
We define the closed subvariety $\La(\alpha, \lambda)$ of $\M(\alpha, \lambda)$ as the fiber $\pi^{-1}(0)$ of the point $0$ under the morphism $\pi$.

\begin{itemize}
\item The variety $\La(\alpha, \lambda)$ is a Lagrangian subvariety of $\M(\alpha, \lambda)$.
In particular, $\dim\La(\alpha, \lambda) = (1/2)\dim\M(\alpha, \lambda)$.

\item The variety $\La(\alpha, \lambda)$ is homotopic to $\M(\alpha, \lambda)$.

\item The variety $\M(\alpha, \lambda)$ is nonempty if and only if so is $\La(\alpha, \lambda)$, and it is equivalent to that the weight space of $V(\lambda)$ of weight $\lambda - \alpha$ is nonzero.

\item If $\M(\alpha, \lambda)$ is nonempty then the dimension of $\M(\alpha, \lambda)$ is equal to $(\lambda, \lambda)-(\lambda - \alpha, \lambda - \alpha)$.

\item For $\alpha, \beta \in Q_+$ with $\alpha \leq \beta$, we have a closed embedding $\M_0(\alpha, \lambda) \hookrightarrow \M_0(\beta, \lambda)$ which sends $0$ to $0$.
Moreover this embedding is an identity for sufficiently large $\alpha, \beta$.
\end{itemize}

We put $d_{\alpha} = \dim \M(\alpha, \lambda)$ for $\alpha \in Q_+$, here $\lambda \in P_+$ is fixed and omitted in the notation.
We put
\[
	\M(\lambda) = \bigsqcup_{\alpha \in Q_+} \M(\alpha, \lambda),
\]
\[
	\La(\lambda) = \bigsqcup_{\alpha \in Q_+} \La(\alpha, \lambda)
\]
and
\[
	\M_0(\lambda) = \bigcup_{\alpha \in Q_+} \M_0(\alpha, \lambda).
\]
Note that the symbol $\M_0(\infty, \mathbf{w})$ was used in \cite{MR1808477} instead of  $\M_0(\lambda)$. 
The composite $\M(\alpha, \lambda) \to \M_0(\alpha, \lambda) \hookrightarrow \M_0(\lambda)$ of the morphisms is also simply denoted by $\pi$.
We have a stratification $\M_0(\lambda) = \bigsqcup_{\alpha \in Q_+} \M_0^{\mathrm{reg}}(\alpha, \lambda)$, where each stratum $\M_0^{\mathrm{reg}}(\alpha, \lambda)$ satisfies the following.

\begin{itemize}
\item For $\lambda \in P_+$ and $\alpha \in Q_+$, $\M_0^{\mathrm{reg}}(\alpha, \lambda)$ is a possibly empty open subset of $\M_0(\alpha, \lambda)$.
It is a nonsingular locally closed subvariety of $\M_0(\lambda)$. 

\item For $\lambda \in P_+$ and $\alpha \in Q_+$, $\M_0^{\mathrm{reg}}(\alpha, \lambda)$ is nonempty if and only if $\lambda - \alpha \in P_+$.

\item The morphism $\pi$ induces an isomorphism between $\pi^{-1}(\M_0^{\mathrm{reg}}(\alpha, \lambda))$ and $\M_0^{\mathrm{reg}}(\alpha, \lambda)$.
If $\M_0^{\mathrm{reg}}(\alpha, \lambda)$ is nonempty then $\pi^{-1}(\M_0^{\mathrm{reg}}(\alpha, \lambda))$ is a dense subset of $\M(\alpha, \lambda)$.
\end{itemize}

\begin{lem}\label{lem:maxdim}
When $\lambda \in P_+$ is fixed and $\alpha \in Q_+$ varies, the dimension of $\M(\alpha, \lambda)$ takes the maximum if and only if $\lambda - \alpha$ is $W$-conjugate to $\lambda_{\mathrm{min}}$.
In particular, the maximum is equal to $(\lambda, \lambda)-(\lambda_{\mathrm{min}}, \lambda_{\mathrm{min}})$.
\end{lem}
\begin{proof}
It is enough to consider the value $(\lambda, \lambda) - (\lambda - \alpha, \lambda - \alpha)$ for various $\alpha \in Q_+$ such that $\lambda - \alpha \in P_+$ since the bilinear form $(\, ,\,)$ is $W$-invariant.
We easily see for $\mu, \nu \in P_+$ that if $\mu < \nu$ then $(\mu, \mu) < (\nu, \nu)$.
Thus by Lemma~\ref{lem:Demazure} (i), $\max\{d_{\alpha} \mid \alpha \in Q_+, \lambda - \alpha \in P_+\}$ is given precisely when $\lambda - \alpha$ is equal to $\lambda_{\mathrm{min}}$.
This completes the proof.
\end{proof}

\begin{lem}\label{lem:dim}
The dimension of $\M_0(\lambda)$ is equal to $\dim \M(\lambda-\lambda_{\mathrm{min}}, \lambda) = (\lambda, \lambda)-(\lambda_{\mathrm{min}}, \lambda_{\mathrm{min}})$.
\end{lem}
\begin{proof}
Since we have the stratification $\M_0(\lambda) = \bigsqcup_{\alpha \in Q_+} \M_0^{\mathrm{reg}}(\alpha, \lambda)$, the dimension of $\M_0(\lambda)$ is equal to $\max\{\dim \M_0^{\mathrm{reg}}(\alpha, \lambda) \mid \alpha \in Q_+\}$.
For $\lambda \in P_+$ and $\alpha \in Q_+$ such that $\M_0^{\mathrm{reg}}(\alpha, \lambda)$ is nonempty, which is equivalent to $\lambda - \alpha \in P_+$, we have
\begin{align*}
	\dim \M_0^{\mathrm{reg}}(\alpha, \lambda) &= \dim \pi^{-1}(\M_0^{\mathrm{reg}}(\alpha, \lambda))\\
						  &= \dim \M(\alpha, \lambda).
\end{align*}
By Lemma~\ref{lem:maxdim}, the maximum of the above is given when $\lambda-\alpha$ is equal to $\lambda_{\mathrm{min}}$.
This completes the proof.
\end{proof}

\subsection{Standard modules}

In this subsection we recall the construction of standard modules for the current Lie algebra $\g[z]$ via the quiver varieties, deduced from \cite{MR1818101} by Varagnolo.
For an algebraic variety $X$ we denote by $H_{\bullet}(X)$ the Borel-Moore homology group of $X$ with complex coefficients.
We simply call it the homology group of $X$ in the sequel.

For $\lambda \in P_+$ and $\alpha, \beta \in Q_+$, we define the variety $Z(\alpha,\beta;\lambda)$ as the fiber product
\[
	Z(\alpha,\beta;\lambda) = \M(\alpha,\lambda) \times_{\M_0(\lambda)} \M(\beta,\lambda)
\]
and $Z(\lambda)$ similarly as
\[
	Z(\lambda) = \M(\lambda) \times_{\M_0(\lambda)} \M(\lambda).
\]
Obviously we have $Z(\lambda) = \bigsqcup_{\alpha, \beta \in Q_+} Z(\alpha,\beta;\lambda)$.
By the general theory of the convolution product (See \cite[2.7]{MR1433132}), the homology group $H_{\bullet}(Z(\lambda))$ has a $\C$-algebra structure via convolution and $H_{\bullet}(\La(\lambda))$ is its module.

For each $\lambda \in P_+$, Varagnolo constructed in \cite[Section~4 Theorem]{MR1818101} an algebra morphism from the Yangian associated with $\g$ to the equivariant homology group of $Z(\lambda)$ with respect to a group action, which is an analogy of a result of Nakajima in \cite{MR1808477} for quantum loop algebras and the equivariant $K$-homology groups.
This makes the equivariant homology group of $\La(\lambda)$ a module over the Yangian.
By forgetting the group action, we deduce the following.

\begin{thm}\label{thm:algebra}
For each $\lambda \in P_+$, there exists a morphism from $U(\g[z])$ to the homology group $H_{\bullet}(Z(\lambda))$ of $\C$-algebras.
In particular, $H_{\bullet}(\La(\lambda))$ is endowed with a $U(\g[z])$-module structure.
\end{thm}

We denote by $M(\lambda)$ the $U(\g[z])$-module $H_{\bullet}(\La(\lambda))$ and call it the \emph{standard module}.
We see that the subspace $H_{\bullet}(\La(\alpha, \lambda))$ corresponds to the $\h$-weight space of $M(\lambda)$ of weight $\lambda - \alpha$ from the explicit definition of the morphism stated in Theorem~\ref{thm:algebra}.

Although the following seems to be known, we give a proof for completeness.

\begin{prop}\label{prop:standard}
We have an isomorphism $M(\lambda) \cong W(\lambda)$ of $U(\g[z])$-modules.
\end{prop}
\begin{proof}
By an argument similar to that in \cite[Proposition~13.3.1]{MR1808477}, we can show that $M(\lambda)$ is generated by the one-dimensional subspace $H_{\bullet}(\La(0, \lambda))$ and that a generator taken from this space satisfies the defining relations of $W(\lambda)$.
Therefore there exists a surjective morphism from $W(\lambda)$ to $M(\lambda)$.
Thus it suffices to show that their dimensions are equal.
The dimension of $M(\lambda)$ is equal to the dimension of $H_{\bullet}(\La(\lambda))$ by the definition and it is known by \cite[Theorem~14.1.2]{MR1808477} to be equal to the product of the dimensions of the fundamental representations of the corresponding quantum loop algebra.
Then by \cite[Corollary~2]{MR2323538}, it is equal to the dimension of $W(\lambda)$.
This completes the proof.
\end{proof}

\section{Cohomological gradings}\label{section:cohomological}

The most part of this section is taken from \cite{MR1433132}, in which a general theory applicable to our setting is developed.
See \cite{MR1433132} for details.

\subsection{Gradings on the convolution algebra and the standard module}

We introduce graded structures into the convolution algebra and its standard module given in Section~\ref{section:quiver}.
Put $A = H_\bullet(Z(\lambda))$ and $M = H_\bullet(\La(\lambda))$.
A formula in \cite[(2.7.9)]{MR1433132} for the change of the degree of the homology groups by the convolution is as follows:
\[
	H_{k}(Z(\alpha, \beta; \lambda)) \times H_{l}(Z(\beta, \gamma; \lambda)) \to H_{k+l-2d_\beta}(Z(\alpha, \gamma; \lambda)),
\]
\[
	H_{k}(Z(\alpha, \beta; \lambda)) \times H_{l}(\La(\beta, \lambda)) \to H_{k+l-2d_\beta}(\La(\alpha, \lambda)).
\]
Put
\[
	H_{[k]}(Z(\lambda)) = \bigoplus_{\alpha, \beta \in Q_+}H_{d_\alpha + d_\beta - k}(Z(\alpha, \beta; \lambda))
\]
and
\[
	H_{[k]}(\La(\lambda)) = \bigoplus_{\alpha \in Q_+}H_{d_\alpha - k}(\La(\alpha, \lambda))
\]
for each $k$.
It is easy to see by the above formula that $A$ is a graded $\C$-algebra by $A = \bigoplus_k H_{[k]}(Z(\lambda))$ and that $M$ is a graded $A$-module by $M = \bigoplus_k H_{[k]}(\La(\lambda))$ \cite[Lemma~8.9.5 and Proposition~8.9.9 (a)]{MR1433132}.

\subsection{Sheaf theoretic analysis of the convolution algebra}

For an algebraic variety $X$ we denote by $D^b(X)$ the bounded derived category of constructible complexes of sheaves on $X$.
The constant sheaf $\C_X$ on $X$ is regarded as an object of $D^b(X)$ concentrated in the degree-zero term.

We abbreviate $\pi_* \colon D^b(\M(\alpha, \lambda)) \to D^b(\M_0(\lambda))$ for the right derived functor of the proper pushforward of $\pi \colon \M(\alpha, \lambda) \to \M_0(\lambda)$.
We denote by $\Ext^{\bullet}(\mathcal{F}, \mathcal{G})$ the $\Ext$ group for objects $\mathcal{F}, \mathcal{G}$ of $D^b(\M_0(\lambda))$.

Put $\mathbb{L}_\alpha = \pi_*\C_{\M(\alpha, \lambda)}[d_\alpha]$ for $\alpha \in Q_+$ and $\mathbb{L} = \bigoplus_{\alpha \in Q_+}\mathbb{L}_\alpha$, which are objects of $D^b(\M_0(\lambda))$.
The $\Ext$ group $\Ext^{\bullet}(\mathbb{L}, \mathbb{L})$ has a graded $\C$-algebra structure by the Yoneda product.
Each $k$ and $\alpha, \beta \in Q_+$, we have an isomorphism $H_{d_\alpha + d_\beta - k}(Z(\alpha, \beta; \lambda)) \cong \Ext^{k}(\mathbb{L}_\alpha, \mathbb{L}_\beta)$ of $\C$-vector spaces \cite[Lemma~8.6.1]{MR1433132}.
Hence $H_{[k]}(Z(\lambda))$ and $\Ext^{k}(\mathbb{L}, \mathbb{L})$ are isomorphic.
In fact, this isomorphism induces that of graded $\C$-algebras \cite[Theorem~8.6.7]{MR1433132}.

\begin{prop}
There exists an isomorphism $H_{\bullet}(Z(\lambda)) \cong \Ext^{\bullet}(\mathbb{L}, \mathbb{L})$ of graded $\C$-algebras.
\end{prop}

Let $i_0 \colon \{0\} \to \M_0(\lambda)$ be the inclusion.
We denote by $i_0^! \colon D^b(\M_0(\lambda)) \to D^b(\{0\})$ the right adjoint of the right derived functor of the proper pushforward of $i_0$.
The cohomology group $H^{\bullet}(i_0^!\mathbb{L})$ is a graded $\Ext^{\bullet}(\mathbb{L}, \mathbb{L})$-module.
We have an isomorphism $H_{d_\alpha - k}(\La(\alpha, \lambda)) \cong H^{k}(i_0^!\mathbb{L_\alpha})$ of $\C$-vector spaces \cite[Lemma~8.5.4]{MR1433132} and this induces an isomorphism $H_{\bullet}(\La(\lambda)) \cong H^{\bullet}(i_0^!\mathbb{L})$ which is compatible with their graded module structures \cite[Proposition~8.6.16]{MR1433132}.
A precise statement is the following.

\begin{prop}
For each $k$, there exists an isomorphism $H_{[k]}(\La(\lambda)) \cong H^{k}(i_0^!\mathbb{L})$ of $\C$-vector spaces such that the following diagram commutes:
\[
	\xymatrix{
	H_{[k]}(Z(\lambda)) \times H_{[l]}(\La(\lambda)) \ar[r]\ar[d]^{\cong} & H_{[k+l]}(\La(\lambda)) \ar[d]^{\cong} \\
	\Ext^k(\mathbb{L}, \mathbb{L}) \times H^l(i_0^!\mathbb{L}) \ar[r] & H^{k+l}(i_0^!\mathbb{L}).
	}
\]
\end{prop}

By the fact that $\pi \colon \M(\alpha, \lambda) \to \M_0(\lambda)$ is \emph{semismall}, proved by Nakajima in \cite[Corollary~10.11]{MR1604167}, $\mathbb{L}$ is a perverse sheaf on $\M_0(\lambda)$.
Then by a property of perverse sheaves, we have $\Ext^k(\mathbb{L}, \mathbb{L}) = 0$ for $k < 0$.
Namely $A = H_{\bullet}(Z(\lambda)) \cong \Ext^{\bullet}(\mathbb{L}, \mathbb{L})$ is positively graded.
By the decomposition theorem together with semismallness of $\pi$, the complex $\mathbb{L}$ decomposes into a direct sum of simple perverse sheaves on $\M_0(\lambda)$ \cite[Theorem~8.4.8 and Proposition~8.9.3]{MR1433132}.
It was proved by Nakajima in \cite[Proposition~15.3.2, see also Theorem~14.3.2 (1)]{MR1808477} that $\mathbb{L}$ decomposes as
\[
	\mathbb{L} \cong \bigoplus_{\substack{\mu \in P_+\\ \mu \leq \lambda}} V_{\mu} \otimes IC(\M_0^{\mathrm{reg}}(\lambda-\mu, \lambda), \C_{\M_0^{\mathrm{reg}}(\lambda-\mu, \lambda)}),
\]
where $IC(\M_0^{\mathrm{reg}}(\lambda-\mu, \lambda), \C_{\M_0^{\mathrm{reg}}(\lambda-\mu, \lambda)})$ is the intersection cohomology complex associated with the constant sheaf on the stratum $\M_0^{\mathrm{reg}}(\lambda-\mu, \lambda)$ and $V_{\mu}$ is a finite-dimensional $\C$-vector space counting its multiplicity.
The decomposition implies that $A_0 = \Ext^0(\mathbb{L}, \mathbb{L}) \cong \bigoplus_{\mu}\End(V_{\mu})^{\oplus m_{\mu}}$ and hence $A_0$ is a finite-dimensional semisimple $\C$-algebra.
Here $m_{\mu}$ denotes the number of the irreducible components of $\M_0^{\mathrm{reg}}(\lambda-\mu, \lambda)$ (In fact we see that $m_{\mu} = 1$ for every $\mu$ since $\M(\lambda-\mu, \lambda)$ is known to be connected.
However connectedness of quiver varieties is not used for our aim.
See Corollary~\ref{cor:connected} and Remark~\ref{rem:connected}).
Each $V_{\mu}$ is regarded as a simple $A$-module through the projection $A \to A_0$.
Then the simple $A$-module $V_{\mu}$ is regarded as a $U(\g[z])$-module through the morphism $U(\g[z]) \to A= H_{\bullet}(Z(\lambda))$ in Theorem~\ref{thm:algebra}.
It is simple as a $U(\g[z])$-module and is isomorphic to $V(\mu)$ by \cite[Theorem~14.3.2 (3)]{MR1808477}.

The following two objects, the $t$-analog of the character and the Kazhdan-Lusztig type polynomial, were introduced by Nakajima in \cite{MR2144973}.

\begin{dfn}
We define the $t$-analog $\chi_t(M(\lambda))$ of the character of the standard module $M(\lambda)$ by
\[
	\chi_t(M(\lambda)) = \sum_{\alpha \in Q_+} \sum_{k=0}^{d_{\alpha}} \dim H_k(\La(\alpha, \lambda))t^{d_{\alpha}-k} e^{\lambda - \alpha}.
\]
\end{dfn}

\begin{dfn}
We define the polynomial $Z_{\lambda \mu}(t)$ by
\[
	Z_{\lambda \mu}(t) = \sum_{k} \dim H^k(i_0^!IC(\M_0^{\mathrm{reg}}(\lambda - \mu, \lambda), \C_{\M_0^{\mathrm{reg}}(\lambda - \mu, \lambda)}))t^k
\]
and call it the Kazhdan-Lusztig type polynomial for the stratum $\M_0^{\mathrm{reg}}(\lambda - \mu, \lambda)$ of the quiver variety.
\end{dfn}

The standard module $M(\lambda)$ has a graded module structure as $M(\lambda) = M = H^{\bullet}(i_0^!\mathbb{L})$ over the positively graded $\C$-algebra $A = \Ext^{\bullet}(\mathbb{L}, \mathbb{L})$.
The associated grading filtration is semisimple since $A_0 = \Ext^{0}(\mathbb{L}, \mathbb{L})$ is semisimple.
By applying the functor $H^{\bullet}(i_0^!(-))$ to $\mathbb{L}$, we see that each coefficient $\dim H^k(i_0^!IC(\M_0^{\mathrm{reg}}(\lambda - \mu, \lambda), \C_{\M_0^{\mathrm{reg}}(\lambda - \mu, \lambda)}))$ of $Z_{\lambda \mu}(t)$ gives the composition multiplicity of $V(\mu)$ in $H^{k}(i_0^!\mathbb{L})$.

It is known by \cite[Theorem~7.3.5]{MR1808477} that $H_{[2k+1]}(\La(\lambda))$ vanishes for every $k$.
We redefine the grading on $M(\lambda)$ by $M(\lambda)_k = H_{[2k]}(\La(\lambda))$, removing the superfluous odd terms.
Although vanishing of $H_{[2k+1]}(Z(\lambda))$ has not been proved so far, this does not affect the module structure of $M(\lambda)$, for every odd term $H_{[2k+1]}(Z(\lambda))$ acts by zero on $M(\lambda)$ as $H_{[2k+1]}(\La(\lambda)) = 0$. 
Under this new grading, the length of the grading filtration on $M(\lambda)$ is given as follows.

\begin{prop}\label{prop:length2}
The length of the grading filtration on the standard module $M(\lambda)$ is equal to $(1/2)\dim \M_0(\lambda) +1= (1/2)\dim \M(\lambda-\lambda_{\mathrm{min}}, \lambda) +1 = (1/2)\{(\lambda, \lambda)-(\lambda_{\mathrm{min}}, \lambda_{\mathrm{min}})\}+1$.
\end{prop}
\begin{proof}
By the definition of the grading, we see that the length of the grading filtration is equal to $(1/2)\max\{d_\alpha \mid \alpha \in Q_+\}+1$.
We have $\max\{d_\alpha \mid \alpha \in Q_+\} = d_{\lambda - \lambda_{\mathrm{min}}} = (\lambda, \lambda)-(\lambda_{\mathrm{min}}, \lambda_{\mathrm{min}})$ by Lemma~\ref{lem:maxdim} and this is equal to $(1/2)\dim \M_0(\lambda)$ by Lemma~\ref{lem:dim}.
\end{proof}

\section{Coincidence of the gradings and its applications}\label{section:coincidence}

\subsection{Coincidence of the gradings}

The following is the second main theorem of this article.

\begin{thm}\label{thm:coincidence}
Under the isomorphism $M(\lambda) \cong W(\lambda)$, the gradings on the both sides coincide.
The associated grading filtration gives a unique Loewy series.
The Loewy length is equal to $(1/2)\dim \M_0(\lambda) +1= (1/2)\dim \M(\lambda-\lambda_{\mathrm{min}}, \lambda) +1 = (1/2)\{(\lambda, \lambda)-(\lambda_{\mathrm{min}}, \lambda_{\mathrm{min}})\}+1$.
\end{thm}
\begin{proof}
We prove that the two gradings coincide.
Recall that $W(\lambda)$ is rigid by Theorem~\ref{thm:rigid} and hence its Loewy series is unique. 
Since the grading filtration on $M(\lambda)$ is semisimple, it suffices to check that its length coincides with the Loewy length of $W(\lambda)$.
This follows from Proposition~\ref{prop:length} (ii) and Proposition~\ref{prop:length2}.
The all remaining assertions are now proved.
\end{proof}

\begin{rem}\label{rem:Varagnolo}
In fact coincidence of the gradings in Theorem~\ref{thm:coincidence} can be deduced also from the fact that the morphism $U(\g[z]) \to H_{\bullet}(Z(\lambda))$ in Theorem~\ref{thm:algebra} is that of graded $\C$-algebras, which is proved by checking directly that the images of the homogeneous generators of $U(\g[z])$ have the appropriate degree.
Although this was not stated in \cite{MR1818101} explicitly, it seems to be known to specialists.
The authors were informed of it by Hiraku Nakajima.
\end{rem}

\begin{cor}\label{cor:connected}
The quiver variety $\M(\alpha, \lambda)$ of type $ADE$ and its Lagrangian subvariety $\La(\alpha, \lambda)$ are connected if they are nonempty.
\end{cor}
\begin{proof}
Let $\alpha$ be an element of $Q_+$ such that $\La(\alpha, \lambda)$ is nonempty. 
Since $\La(\alpha, \lambda)$ is homotopic to $\M(\alpha, \lambda)$, it suffices to show that $\dim H_0(\La(\alpha, \lambda)) = 1$.
Recall that $W(\lambda) \cong L_w(\Lambda)$ where $w_0w\Lambda = \lambda + \Lambda_0$.
By Theorem~\ref{thm:coincidence}, we see that $H_0(\La(\alpha, \lambda))$ corresponds to the $\hat{\h}$-weight space of $L_w(\Lambda)$ of weight $(\lambda - \alpha) + \Lambda_0 + (d_{\alpha}/2)\delta$ under the isomorphism $M(\lambda) \cong W(\lambda) \cong L_w(\Lambda)$.
This weight is equal to $t_{-\alpha}(\lambda + \Lambda_0) = t_{-\alpha}w_0w\Lambda$, an extremal weight of $L(\Lambda)$, and hence its weight space is one-dimensional. 
This completes the proof.
\end{proof}

\begin{rem}\label{rem:connected}
Connectedness of quiver varieties for arbitrary graphs has been already proved by Crawley-Boevey in \cite{MR1834739}.
Here we deduce this fact for the quiver varieties of type $ADE$ as a corollary of our result.
\end{rem}

\subsection{The gradings on Weyl modules and crystals}\label{subsection:crystal}

The grading on the standard module $M(\lambda)$ comes from the degree of the homology group of the quiver variety as explained in Section~\ref{section:cohomological}.
The grading on the Weyl module $W(\lambda)$ is related to the energy function defined on a certain crystal, which was proved by the second author in \cite{naoi}.
In this subsection, we review this subject.
Together with the fact that the both gradings coincide as proved in Theorem~\ref{thm:coincidence}, we obtain nontrivial relations between quiver varieties and crystals, Corollary~\ref{cor:Poincare} and Corollary~\ref{cor:KL}.

Let $M$ be a finite-dimensional graded $U(\g[z])$-module.
We denote by $M_{\mu}$ the $\h$-weight space of $M$ of weight $\mu \in P$.
Since the action of $\h$ is degree zero, we have $M=\bigoplus_{k \in \z}\bigoplus_{\mu \in P} (M_k \cap M_{\mu})$.
The \emph{graded character} $\gch_t M$ of $M$ is defined by
\[
	\gch_t M = \sum_{k \in \z}\sum_{\mu \in P} \dim(M_k \cap M_{\mu})t^k e^{\mu}.
\]
This is a refinement of the usual $\h$-character of $M$ as $\gch_1 M = \ch M$.
The \emph{graded composition multiplicity} of the simple module $V(\mu)$ in $M$ is the polynomial defined by
\[
	\sum_{k \in \z}[M_k:V(\mu)]t^k,
\]
where $[M_k:V(\mu)]$ denotes the composition multiplicity of $V(\mu)$ in $M_k$.
We remark that the $t$-analog $\chi_t(M(\lambda))$ of the character and the Kazhdan-Lusztig type polynomial $Z_{\lambda \mu}(t)$ introduced in Section~\ref{section:cohomological} are nothing but the graded character of $M(\lambda)$ and the graded composition multiplicity of $V(\mu)$ in $M(\lambda)$ respectively, where the grading is original one.
Hence we have
\[
	\chi_t(M(\lambda)) = \gch_{t^2} W(\lambda)
\]
and 
\[
	Z_{\lambda \mu}(t) = \sum_{k \geq 0}[W(\lambda)_k:V(\mu)]t^{2k}
\]
by coincidence of the gradings in Theorem~\ref{thm:coincidence}.

Let us turn to crystals.
We refer to the seminal paper by Kashiwara \cite{MR1115118} and a standard text book \cite{MR1881971} for the basic theory of crystals and to \cite{MR1890649} for the level-zero fundamental representations of the quantum affine algebra $U_q'(\hat{\g})$ and their crystal bases.
In the sequel we treat only $\hat{P}_{\cl}$-crystals, where $\hat{P}_{\mathrm{cl}}=\hat{P}/\z\delta$ is the weight lattice of $U_q'(\hat{\g})$.
Note that $P$ can be regarded as a subset of $\hat{P}_{\mathrm{cl}}$.
Then for a $\hat{P}_{\cl}$-crystal with level-zero weight, the weight map $\wt$ takes its values in $P$.

The definition of one-dimensional sums was given by Hatayama, Kuniba, Okado, Takagi and Yamada in \cite{MR1745263} and by Hatayama, Kuniba, Okado, Takagi and Tsuboi in \cite{MR1903978} motivated by the study of solvable lattice models.
Here we need them only for tensor products of the crystal bases of level-zero fundamental representations, while they are defined for a wider class of crystals in general.
Our treatment follows an approach by Naito and Sagaki in \cite{MR2474320}.
We denote by $B^i$ the crystal base of the level-zero fundamental representation of the quantum affine algebra $U_q'(\hat{\g})$ associated with the fundamental weight $\varpi_i$ for $i \in I$. 
For a given sequence $\mathbf{i} = (i_1, \ldots, i_l)$ of elements of $I$ we put $B_{\mathbf{i}} = B^{i_1} \otimes \cdots \otimes B^{i_l}$.
Let $D_{\bfi}$ be the \emph{energy function} defined on the crystal $B_{\bfi}$ and $\dext$ the extra constant, which is a certain integer, as in \cite[Subsection~4.1]{MR2474320}.
We denote by $\mathbb{B}(\lambda)_{\cl}$ the crystal consisting of all Lakshmibai-Seshadri paths of shape $\lambda \in P_+$ modulo $\delta$, which was proved by Naito and Sagaki in \cite[Corollary~4.4]{MR2146858} to be isomorphic to $B_{\bfi}$ for every $\bfi$ satisfying $\sum_{k=1}^{l}\varpi_{i_k} = \lambda$.
Let $D$ be the degree function on $\mathbb{B}(\lambda)_{\cl}$ defined in \cite[Subsection~3.1]{MR2474320}.
The following theorem was proved in \cite[Thorem~4.1.1]{MR2474320}.

\begin{thm}\label{thm:degree}
Under the isomorphism $B_{\bfi} \cong \mathbb{B}(\lambda)_{\cl}$, the function $D_{\bfi} - \dext$ on $B_{\bfi}$ corresponds to the degree function $D$ on $\mathbb{B}(\lambda)_{\cl}$.
\end{thm}

Following \cite{MR1745263}, \cite{MR1903978} and \cite{MR2474320} we introduce a polynomial called the one-dimensional sum and its variant.

\begin{dfn}
Let $\mathbf{i}$ be a sequence of elements of $I$.
We define the \emph{one-dimensional sum} $X(B_{\bfi}, \mu, t)$ associated with the crystal $B_{\bfi}$ and $\mu \in P_+$ by
\[
X(B_{\bfi}, \mu, t) = \sum_{\substack{b \in B_{\bfi}\\ \e{i}b=0\, (i \in I)\\ \wt b = \mu}} t^{D_{\bfi}(b)}.
\]
For $\lambda, \mu \in P_+$, we define the polynomial $X(\lambda, \mu, t)$ by
\[
X(\lambda, \mu, t) = \sum_{\substack{b \in \mathbb{B}(\lambda)_{\mathrm{cl}}\\ \e{i}b=0\, (i \in I)\\ \wt b = \mu}} t^{D(b)}.
\]
\end{dfn}

\begin{rem}
Let $\lambda$ and $\mu$ be elements of $P_+$ and take a sequence $\mathbf{i} = (i_1, \ldots, i_l)$ of elements of $I$ so that $\sum_{k=1}^{l}\varpi_{i_k} = \lambda$.
Then we have
\[
X(\lambda, \mu, t) = t^{-\dext}X(B_{\bfi}, \mu, t)
\]
by Theorem~\ref{thm:degree}.
\end{rem}

\begin{rem}
In general, $X(\lambda, \mu, t)$ is a polynomial in $t^{-1}$ with its coefficients in nonnegative integers since the degree function takes its values in $\z_{\leq 0}$ by its definition.
\end{rem}

The following theorem was proved by the second author in \cite[Theorem~9.2 and Corollary~9.6]{naoi}, generalizing a result in \cite{MR2271991} which concerns the case of type $A$.

\begin{thm}\label{thm:gch}
Let $\lambda$ be an element of $P_+$.
\begin{enumerate}
\item We have
\begin{align*}
\gch_{t} W(\lambda) &= \sum_{b \in \mathbb{B}(\lambda)_{\cl}} t^{-D(b)} e^{\wt b}\\
&= \sum_{\mu \in P_+} X(\lambda, \mu, t^{-1})\ch V(\mu).
\end{align*}
\item For $\mu \in P_+$, we have
\[
	\sum_{k \geq 0}[W(\lambda)_k:V(\mu)]t^k = X(\lambda, \mu, t^{-1}).
\]
\end{enumerate}
\end{thm}

We obtain the following corollary of Theorem~\ref{thm:coincidence} and Theorem~\ref{thm:gch}, a formula for the Poincar\'{e} polynomials of  quiver varieties, as an equality for the graded character.

\begin{cor}\label{cor:Poincare}
We have
\begin{align*}
\sum_{\alpha \in Q_+} \sum_{k=0}^{d_{\alpha}} \dim H_k(\La(\alpha, \lambda))t^{d_{\alpha}-k} e^{\lambda - \alpha} &= \sum_{b \in \mathbb{B}(\lambda)_{\cl}}t^{-2D(b)}e^{\wt b}\\
&= \sum_{\mu \in P_+} X(\lambda, \mu, t^{-2})\ch V(\mu).
\end{align*}
In particular, we have
\[
\sum_{k=0}^{d_{\alpha}} \dim H_k(\La(\alpha, \lambda))t^{d_{\alpha}-k}= \sum_{\substack{b \in \mathbb{B}(\lambda)_{\cl}\\ \wt b = \lambda - \alpha}}t^{-2D(b)}.
\]
\end{cor}

We also obtain an equality for the graded composition multiplicity.
Recall that $Z_{\lambda \mu}(t)$ denotes the Kazhdan-Lusztig type polynomial for the stratum $\M_0^{\mathrm{reg}}(\lambda-\mu, \lambda)$ of the quiver variety.

\begin{cor}\label{cor:KL}
We have
\[
	Z_{\lambda \mu}(t) = X(\lambda, \mu, t^{-2}).
\]
\end{cor}

\subsection{Lusztig's fermionic conjecture and the $X=M$ conjecture}

Lusztig conjectured in \cite{lus} that the Poincar\'{e} polynomials of quiver varieties are described in terms of the fermionic forms also introduced in \cite{MR1745263} and \cite{MR1903978}.
The so-called $X=M$ conjecture asserts that the one-dimensional sum and the fermionic form coincide up to some constant power.
We discuss these subjects in this subsection.

For a sequence $\mathbf{m}=(m_k^{(i)})_{i \in I, k \in \z_{>0}}$ of nonnegative integers with finitely many nonzero terms and $\lambda \in P_+$, we define the integer $p_k^{(i)}(\mathbf{m}, \lambda)$ for $i \in I$ and $k \in \z_{>0}$ by
\[
	p_k^{(i)}(\mathbf{m}, \lambda) = \langle h_i, \lambda \rangle - \sum_{j \in I}(\alpha_i, \alpha_j)\sum_{l \geq 1} \min\{k,l\}m_l^{(j)}
\]
and the integer $c(\mathbf{m}, \lambda)$ by 
\[
	c(\mathbf{m}, \lambda) = \frac{1}{2}\sum_{i,j \in I} (\alpha_i, \alpha_j)\sum_{k,l \geq 1} \min\{k,l\}m_k^{(i)}m_l^{(j)} - \sum_{i \in I}\sum_{k \geq 1} \langle h_i, \lambda \rangle m_k^{(i)}.
\]
For $\lambda, \mu \in P_+$ we define the subset $S(\lambda, \mu)$ of $\z_{\geq 0}^{(I \times \z_{>0})}$ by
\[
	S(\lambda, \mu) = \biggl\{\mathbf{m} \biggm| \sum_{i \in I}\sum_{k \geq 1} km_k^{(i)}\alpha_i = \lambda-\mu \biggr\}.
\]
With the notation as above, we define the fermionic form $M(\lambda, \mu, t)$ associated with $\lambda, \mu \in P_+$ by
\[
	M(\lambda, \mu, t) = \sum_{\mathbf{m} \in S(\lambda, \mu)}t^{c(\mathbf{m},\lambda)}\prod_{i \in I, k\geq1} \gauss{p_k^{(i)}(\mathbf{m},\lambda) + m_k^{(i)}}{m_k^{(i)}}_t,
\]
where the Gaussian binomial coefficients are defined in a usual manner with $t$-integers as $[n]_t=(t^n - 1)/(t-1)$ and we set
\[
	\gauss{m}{n}_{t}=0
\] for $m<n$.
Then $M(\lambda, \mu, t)$ turns out to be a polynomial in $t^{-1}$ with its coefficients in nonnegative integers.
Lusztig's fermionic conjecture \cite[Conjecture~A]{lus} is the following:
\[
	\sum_{\alpha \in Q_+} \sum_{k=0}^{d_{\alpha}} \dim H_k(\La(\alpha, \lambda))t^{d_{\alpha}-k} e^{\lambda - \alpha} = \sum_{\mu \in P_+} M(\lambda, \mu, t^{-2})\ch V(\mu).
\]
As pointed out in \cite[Section~9]{naoi}, results of Ardonne and Kedem in \cite{MR2290922} and Di Francesco and Kedem in \cite{MR2428305}, together with \cite{MR2323538}, imply that the graded composition multiplicity of the simple module $V(\mu)$ in the Weyl module $W(\lambda)$ coincides with the fermionic form $M(\lambda, \mu, t^{-1})$.
Then by Theorem~\ref{thm:coincidence}, we conclude that Lusztig's conjecture is true.

Mozgovoy proved in \cite{moz} a variant of the above formula.
It implies with an easy argument the following:
\[
	\sum_{\alpha \in Q_+} \sum_{k=0}^{d_{\alpha}} \dim H_k(\La(\alpha, \lambda))t^{d_{\alpha}-k} e^{\lambda - \alpha} = \sum_{\mu \in P_+} N(\lambda, \mu, t^{-2})\ch V(\mu),
\]
where the polynomial $N(\lambda, \mu, t)$ is defined in a way similar to $M(\lambda, \mu, t)$ with the different binomial coefficients.
It was conjectured in \cite{MR1745263} and \cite{MR1903978} that $M(\lambda, \mu, t) = N(\lambda, \mu, t)$.
The above formulas imply that the conjecture is also true.

Also as explained in \cite[Section~9]{naoi}, the equality $X(\lambda, \mu, t) = M(\lambda, \mu, t)$ holds by combining the results mentioned above and Theorem~\ref{thm:gch}.
Indeed both $X(\lambda,\mu,t^{-1})$ and $M(\lambda,\mu,t^{-1})$ coincide with the graded composition multiplicity of $V(\mu)$ in $W(\lambda)$.
Thus the $X=M$ conjecture for the case of tensor products of level-zero fundamental representations has been now solved.

So far we have achieved the equalities:
\[
	Z_{\lambda\mu}(t) = X(\lambda,\mu,t^{-2}) = M(\lambda,\mu,t^{-2}) = N(\lambda,\mu,t^{-2}).
\]
\begin{rem}
One-dimensional sums and fermionic forms are both defined for tensor products of Kirillov-Reshetikhin modules for quantum affine algebras and it was conjectured in \cite{MR1745263} and \cite{MR1903978} that they coincide up to some constant power.
In a general setting, this conjecture is still open.
The authors do not know whether quiver varieties are related to them for general Kirillov-Reshetikhin modules or not.
\end{rem}

\subsection{Type $A$}

In this subsection, we pick up two previous works for type $A$, \cite{MR1285530} by Nakajima and \cite{MR2271991} by Chari and Loktev.
Let $\g$ be the simple Lie algebra of type $A_n$.
The index set $I$ is identified with $\{1, \ldots, n\}$ with a usual numbering as $\langle h_i, \alpha_j \rangle = 2\delta_{i,j} -\delta_{i,j-1}-\delta_{i,j+1}$.

Nakajima calculated in \cite{MR1285530} the Poincar\'{e} polynomials of quiver varieties of type $A$ in terms of tableaux.
We identify the set $P_+$ of all dominant integral weights with the set of all partitions of length less than or equal to $n$ as usual.
Let $T(\alpha, \lambda)$ be the set of all row-increasing tableaux of shape ${}^t\lambda$ and weight $\lambda-\alpha$.
Nakajima defined a certain function $l$ on $T(\alpha, \lambda)$ and proved the following \cite[Theorem~5.15]{MR1285530}:
\[
	\sum_{k=0}^{d_{\alpha}} \dim H_k(\La(\alpha, \lambda))t^k = \sum_{T \in T(\alpha, \lambda)} t^{2l(T)}.
\]
By comparing this formula with ours in Corollary~\ref{cor:Poincare}, we see that there exists an abstract bijection between the set $\{b \in \mathbb{B}(\lambda)_{\cl} \mid \wt b = \lambda - \alpha\}$ and $T(\alpha, \lambda)$
 such that the degree function $D$ corresponds to the function $-(1/2)d_{\alpha} + l$.
We should give an explicit description of the bijection.
It is known that $\mathbb{B}(\lambda)_{\cl}$ is identified with the set of all column-increasing tableaux of shape $\lambda$.
By transposing the tableaux, we obtain a bijection.
However this naive one does not satisfy the property mentioned above.
Thus we need further study to understand a relation between his result and ours.

Chari and Loktev gave in \cite{MR2271991} a quite explicit description of the graded characters of Weyl modules for type $A$.
In our viewpoint, it gives an explicit formula for the Poincar\'{e} polynomials of quiver varieties of type $A$. 
Put $\alpha_{i,j}=\sum_{s=i}^{j}\alpha_s$ for $i,j \in I$ with $i \leq j$.
These elements form the positive roots of the root system of type $A_n$.
For $\alpha \in Q_+$ we define the subset $S(\alpha)$ of $\z_{\geq 0}^{n(n+1)/2}$ by
\[
	S(\alpha) = \biggl\{(l_{i,j})_{1 \leq i \leq j \leq n}  \biggm| \sum_{1 \leq i \leq j \leq n} l_{i,j} \alpha_{i,j} = \alpha \biggr\}.
	\]
A formula in \cite[2.1.4 Proposition]{MR2271991} yields the following:
\begin{align*}
&\sum_{k=0}^{d_{\alpha}} \dim H_k(\La(\alpha, \lambda))t^{d_{\alpha}-k}\\
&=\sum_{(l_{i,j}) \in S(\alpha)}\prod_{1 \leq i \leq j \leq n}\gauss{\langle h_i, \lambda \rangle + \sum_{s=j+1}^{n} l_{i+1,s} - \sum_{s=j+1}^{n} l_{i,s}}{l_{i,j}}_{t^2}.
\end{align*}


\begin{thebibliography}{HKOTY}

\bibitem[AK]{MR2290922}
Eddy Ardonne and Rinat Kedem, \emph{Fusion products of {K}irillov-{R}eshetikhin modules and fermionic multiplicity formulas}, J. Algebra \textbf{308} (2007), no.~1, 270--294.

\bibitem[BGS]{MR1322847}
Alexander Beilinson, Victor Ginzburg, and Wolfgang Soergel, \emph{Koszul duality patterns in representation theory}, J. Amer. Math. Soc. \textbf{9} (1996), no.~2, 473--527.

\bibitem[CL]{MR2271991}
Vyjayanthi Chari and Sergei Loktev, \emph{Weyl, {D}emazure and fusion modules for the current algebra of {$\mathfrak{sl}_{r+1}$}}, Adv. Math. \textbf{207} (2006), no.~2, 928--960.

\bibitem[CM]{MR2078944}
Vyjayanthi Chari and Adriano A. Moura, \emph{Spectral characters of finite-dimensional representations of affine algebras}, J. Algebra \textbf{279} (2004), no.~2, 820--839.

\bibitem[CP]{MR1850556}
Vyjayanthi Chari and Andrew Pressley, \emph{Weyl modules for classical and quantum affine algebras}, Represent. Theory \textbf{5} (2001), 191--223 (electronic).

\bibitem[CG]{MR1433132}
Neil Chriss and Victor Ginzburg, \emph{Representation theory and complex geometry}, Birkh\"auser Boston Inc., Boston, MA, 1997.

\bibitem[CB]{MR1834739}
William Crawley-Boevey, \emph{Geometry of the moment map for representations of quivers}, Compos. Math. \textbf{126} (2001), no.~3, 257--293.

\bibitem[DFK]{MR2428305}
Philippe Di~Francesco and Rinat Kedem, \emph{Proof of the combinatorial {K}irillov-{R}eshetikhin conjecture}, Int. Math. Res. Not. (2008), no.~7, Art.~ID rnn006, 57pp.

\bibitem[FeL]{MR2102326}
Boris Feigin and Sergei Loktev, \emph{Multi-dimensional {W}eyl modules and symmetric functions}, Comm. Math. Phys. \textbf{251} (2004), no.~3, 427--445.

\bibitem[FoL]{MR2323538}
Ghislain Fourier and Peter Littelmann, \emph{Weyl modules, {D}emazure modules, {KR}-modules, crystals, fusion products and limit constructions}, Adv. Math. \textbf{211} (2007), no.~2, 566--593.

\bibitem[HKOTT]{MR1903978}
Goro Hatayama, Atsuo Kuniba, Masato Okado, Taichiro Takagi, and Zengo Tsuboi, \emph{Paths, crystals and fermionic formulae}, Math{P}hys odyssey, 2001, Prog. Math. Phys., vol.~23, Birkh\"auser Boston, Boston, MA, 2002, 205--272.

\bibitem[HKOTY]{MR1745263}
Goro Hatayama, Atsuo Kuniba, Masato Okado, Taichiro Takagi, and Yasuhiko Yamada, \emph{Remarks on fermionic formula}, Recent developments in quantum affine algebras and related topics ({R}aleigh, {NC}, 1998), Contemp. Math., vol.~248, Amer. Math. Soc., Providence, RI, 1999, 243--291.

\bibitem[H]{MR2651380}
Tam{\'a}s Hausel, \emph{Kac's conjecture from {N}akajima quiver varieties}, Invent. Math. \textbf{181} (2010), no.~1, 21--37.

\bibitem[HK]{MR1881971}
Jin Hong and Seok-Jin Kang, \emph{Introduction to quantum groups and crystal bases}, Graduate Studies in Mathematics, vol.~42, American Mathematical Society, Providence, RI, 2002.

\bibitem[Kac]{MR1104219}
Victor~G. Kac, \emph{Infinite dimensional {L}ie algebras}, third ed., Cambridge University Press, Cambridge, 1990.

\bibitem[Kas1]{MR1115118}
Masaki Kashiwara, \emph{On crystal bases of the {$Q$}-analogue of universal enveloping algebras}, Duke Math. J. \textbf{63} (1991), no.~2, 465--516.

\bibitem[Kas2]{MR1890649}
\bysame, \emph{On level-zero representations of quantized affine algebras}, Duke Math. J. \textbf{112} (2002), no.~1, 117--175.

\bibitem[Ko]{MR2657446}
Ryosuke Kodera, \emph{Extensions between finite-dimensional simple modules over a generalized current {L}ie algebra}, Transform. Groups \textbf{15} (2010), no.~2, 371--388.

\bibitem[L]{lus}
George Lusztig, \emph{Fermionic form and {B}etti numbers}, arXiv:0005010.

\bibitem[M]{moz}
Sergey Mozgovoy, \emph{Fermionic forms and quiver varieties}, arXiv:0610084.

\bibitem[NS1]{MR2146858}
Satoshi Naito and Daisuke Sagaki, \emph{Crystal of {L}akshmibai-{S}eshadri paths associated to an integral weight of level zero for an affine {L}ie algebra}, Int. Math. Res. Not. (2005), no.~14, 815--840.

\bibitem[NS2]{MR2474320}
\bysame, \emph{Lakshmibai-{S}eshadri paths of level-zero shape and one-dimensional sums associated to level-zero fundamental representations}, Compos. Math. \textbf{144} (2008), no.~6, 1525--1556.

\bibitem[Nak1]{MR1285530}
Hiraku Nakajima, \emph{Homology of moduli spaces of instantons on {ALE} spaces. {I}}, J. Differential Geom. \textbf{40} (1994), no.~1, 105--127.

\bibitem[Nak2]{MR1302318}
\bysame, \emph{Instantons on {ALE} spaces, quiver varieties, and {K}ac-{M}oody algebras}, Duke Math. J. \textbf{76} (1994), no.~2, 365--416.

\bibitem[Nak3]{MR1604167}
\bysame, \emph{Quiver varieties and {K}ac-{M}oody algebras}, Duke Math. J. \textbf{91} (1998), no.~3, 515--560.

\bibitem[Nak4]{MR1808477}
\bysame, \emph{Quiver varieties and finite-dimensional representations of quantum affine algebras}, J. Amer. Math. Soc. \textbf{14} (2001), no.~1, 145--238.

\bibitem[Nak5]{MR2144973}
\bysame, \emph{Quiver varieties and {$t$}-analogs of {$q$}-characters of quantum affine algebras}, Ann. of Math. (2) \textbf{160} (2004), no.~3, 1057--1097.

\bibitem[Nao]{naoi}
Katsuyuki Naoi, \emph{Weyl modules, {D}emazure modules and finite crystals for non-simply laced type}, arXiv:1012.5480.

\bibitem[V]{MR1818101}
Michela Varagnolo, \emph{Quiver varieties and {Y}angians}, Lett. Math. Phys. \textbf{53} (2000), no.~4, 273--283.

\end{thebibliography}
\end{document}